\documentclass{naw-latex}
\usepackage{amssymb}
\usepackage{xcolor}
\newtheorem{theorem}{Theorem}

\newtheorem{conjecture}{Conjecture}

\bibitem{bAlLint}{
P. J. van Albada and J. H. van Lint,  Reciprocal bases for the integers, \textit{Amer. Math. Monthly}, 70 (1963), 170–174.}

\bibitem{Ave}{G. Averkov, Difference Between Families of Weakly and Strongly Maximal Integral Lattice-Free Polytopes, in \emph{Interactions with Lattice Polytopes} (Springer, ILP 2017) Vol. 386., pages 1--10}

\bibitem{BFM}{A. Beltr\'{a}n, M. J. Felipe, and C. Melchor, Landau's theorem on conjugacy classes for normal subgroups, \emph{International Journal of Algebra and Computation} 26 (2016), 1453-1466.}

\bibitem{Blo2021}{T. F. Bloom, On a density conjecture about unit fractions, arXiv:2112.03726 (2021).}

\bibitem{BL}{M. Bright and D. Loughran, Brauer-Manin obstruction for Erd\H{o}s-Straus surfaces. \emph{Bull. Lond. Math. Soc.} 52 (2020), no. 4, 746–761.}

\bibitem{Bowning-Elsholtz}{T. D. Browning and C. Elsholtz. The number of representations of rationals
as a sum of unit fractions. {\emph{Illinois J. Math.}} 55 (2011), no. 2,
685--696.}

\bibitem{Cro2}{E. S. Croot, 
On unit fractions with denominators in short intervals,
\emph{Acta Arith.} 99 (2001), no. 2, 99–114.}

\bibitem{Cro2003}{E. S. Croot, On a coloring conjecture about unit fractions. \textit{Ann. of Math.} (2) 157 (2003), 545-556.}

\bibitem{Dedekind:1931}{R.~Dedekind.
{\em {Gesammelte mathematische Werke, Band 2, Hrsg. v. Robert Fricke,
  Emmy Noether u. Oeystein Ore., \"Uber Zerlegungen von Zahlen durch ihren
  gr\"o\ss ten gemeinsamen Teiler, (Festschrift der Universit\"{a}t
  Braunschweig, 1897)}}. Braunschweig: Friedr. Vieweg \& Sohn A.-G., 1931.
}

\bibitem{Elsholtz:2001}{C. Elsholtz, Sums of $k$ unit fractions,
\emph{Trans. Amer. Math. Soc.} 353 (2001), 3209-3227.} 

\bibitem{Elsholtz:1998}{C. Elsholtz, Sums of $k$ unit fractions, Ph.D. Thesis  (Darmstadt), 1998.
Shaker Verlag, Aachen.}

\bibitem{Elsholtz:2015}{
C. Elsholtz, 
Egyptian fractions with odd denominators, \textit{Q. J. Math.} 67 (2016), no. 3, 425-430.} 

\bibitem{Elsholtz-Heuberger-Krenn}{C.~Elsholtz, C.~Heuberger, and D.~Krenn,
Algorithmic counting of nonequivalent compact Huffman codes,
https://arxiv.org/abs/1901.11343.}

\bibitem{Elsholtz-Heuberger-Prodinger}{C. Elsholtz, C. Heuberger, and H. Prodinger. The number of Huffman
codes, compact trees, and sums of unit fractions. IEEE Trans. Inform.
Theory 59 (2013), no. 2, 1065--1075.}

\bibitem{Elsholtz-Planitzer:2020}{C.~Elsholtz and S.~Planitzer, The number of solutions of the Erd\H{o}s-Straus Equation and sums of k unit fractions, \textit{Proc. R. Soc. Edinb. A: Math.} 150(3) (2020), 1401-1427.}

\bibitem{elsholtz-planitzer:2021}{C.~Elsholtz and S.~Planitzer, Sums of four and more unit fractions and approximate parametrizations, \textit{Bull. Lond. Math. Soc.} 53 (3), 2021, 695-709.}

\bibitem{elsholtz-tao}{C.~Elsholtz and T.~Tao. Counting the number of solutions to the  Erd{\H o}s-{S}traus  equation on unit fractions, \textit{J. Aust. Math. Soc.} 94 (2013), no. 1, 50-105.  }

\bibitem{Eppstein}{D. Eppstein, Ten algorithms for Egyptian fractions,\textit{ Mathematica in Education and Research} 4(2) (1995), 5-15.}

\bibitem{Erdos:1950}{P.~Erd\H{o}s,
{Az ${1}/{x_1} + {1}/{x_2} + \ldots + {1}/{x_n}  ={a}/{b}$ egyenlet eg\'{e}sz sz\'{a}m\'{u} megold\'{a}sair\'{o}l}, \emph{Mat. Lapok} 1 (1950), 192-210.}

\bibitem{ErdosandGraham}{
P. Erd\H{o}s and R. L. Graham,
Old and new problems and results in combinatorial number theory,
\textit{Monographies de L'Enseignement Math\'{e}matique}, 28. Universit\'{e} de Gen\`{e}ve, Geneva, 1980.}

\bibitem{Graham:2013}{R.L.~Graham,
Paul Erd\H{o}s and Egyptian fractions, in \textit{Erdös centennial}, 289-309, Bolyai Soc. Math. Stud., 25, J\'{a}nos Bolyai Math. Soc., Budapest, 2013. }

\bibitem{Guy}{R. K. Guy, Unsolved Problems in Number Theory, Springer, Third Edition, 2004.}

\bibitem{Hipler:2002}{Andreas Hipler,
Zu Modulspannweiten
und zu einem Problem
von Erd\H{o}s und Graham, 
Ph.D. Thesis, 2002, University of Mainz.}

\bibitem{Izhboldin}{O. T. Izhboldin and L. D. Kurlyandchik, Unit fractions, \textit{Proc. St. Petersburg Math. Soc.}, 111 (1995) 193-200.}

\bibitem{Kon2014}{S. V. Konyagin, Double Exponential Lower Bound for the Number of Representations of Unity by Egyptian Fractions, \textit{Math. Notes} 95 (2014), no. 1-2, 277-281.}

\bibitem{Lan1903}{E. Landau, \"Uber die Klassenzahl der bin\"aren quadratischen Formen von negativer Diskriminante, \textit{Math. Ann.} 56 (1903), 671-676.}

\bibitem{Mar}{G. Martin, Denser Egyptian fractions, \textit{ Acta Arith.} 95 (2000), no. 3, 231–260.}

\bibitem{Mor1969}{L. J. Mordell, Diophantine Equations, volume 30 of Pure and Applied Mathematics. Academic Press, 1969.}

\bibitem{Nak1939}{M. Nakayama, On the Decomposition of a Rational Number into ``Stammbr\"{u}che'', \textit{Tohuko J. Math.} 46 (1939), 1-21.}

 \bibitem{NZ11}{B. Nill and G. M. Ziegler, Projecting lattice polytopes without interior lattice points, \textit{Math. Oper. Res.} 36 (2011), no. 3, 462–467.}

\bibitem{Obl1950}{M.R. Obl\'{a}th, {Sur l' \'{e}quation diophantienne ${4}/{n}={1}/{x_1}+{1}/{x_2} +{1}/{x_3}$}, \emph{Mathesis} \textbf{59} (1950), 308-316.}

\bibitem{Rhind}{
The Papyrus Rhind, approximately 1650–1550 BCE, written by the scribe by Ahmes.}

\bibitem{Ros1954}{
L. A. Rosati, Sull'equazione diofantea $\frac{4}{n}=\frac{1}{x_1}+\frac{1}{x_2}+\frac{1}{x_3}$, \textit{Bollettino della Unione Matematica Italiana} (3), 9: 59-64 (1954).
}

\bibitem{Sal}{S. Salez, The Erd\H{o}s-Straus conjecture: New modular equations and checking up to $N=10^{17}$, arXiv:1406.6307.} 

\bibitem{Sandor:2003}{C. S\'{a}ndor,
 On the number of solutions of the Diophantine equation $\sum_{i=1}^n\frac{1}{x_i}=1$,
 \textit{Period. Math. Hungar.} 47 (2003), no. 1-2, 215-219.}

\bibitem{Schinzel:1956}{A. Schinzel, Sur quelques propri\'{e}t\'{e}s des nombres $3/n$ et $4/n$,  o\`{u} $n$ est un nombre impair. \emph{Mathesis} 65 (1956), 219-222.}

\bibitem{Schinzel}{
A. Schinzel, 
Erd\H{o}s's work on finite sums of unit fractions, in \textit{Paul Erd\H{o}s and his mathematics}, I (Budapest, 1999), 629-636,
Bolyai Soc. Math. Stud., 11, J\'{a}nos Bolyai Math. Soc., Budapest, 2002.}

\bibitem{Schinzel:2000}{A. Schinzel,
On sums of three unit fractions with polynomial denominators.
\textit{Funct. Approx. Comment. Math.} 28 (2000), 187-194.}

\bibitem{Shen}{Z. Shen,
On the Diophantine equation $\sum^k_{i=0}1/x_i=a/n$,
\textit{Chinese Ann. Math. Ser. B} 7 (1986), 2, 213-220.}

\bibitem{Sierpinski:1956}{W. Sierpi\'{n}ski, Sur les d\'{e}compositions de nombres rationnels en fractions primaires, \textit{Mathesis} 65 (1956), 16-32.}

\bibitem{Stewart:1964}{B. M. Stewart, Theory of numbers,
(second edition). The Macmillan Company, New York; Collier Macmillan Ltd., London 1964.}

\bibitem{Takenouchi}{T. Takenouchi, On an indeterminate equation, \textit{Proc. Phys.-Math. Soc. Japan} (3), 3 (1921) 78-92.}

\bibitem{Ter1971}{D.G. Terzi, On a conjecture by {Erd\H{o}s-Straus},
{\em Nordisk Tidskr. Informations-Behandling (BIT)} \textbf{11} (1971), 212-216.}
 
\bibitem{Vau1970}{R. Vaughan, {On a problem of Erd\H{o}s, Straus and Schinzel}, \emph{Mathematika} 17 (1970), 193-198.}

\bibitem{Viola}{
C.~Viola,
On the diophantine equations $\Pi ^{k}_{0}x_{i}-\sum
   ^{k}_{0}\,x_{i}=n$ and $\sum ^{k}_{0}\,1/x_{i}=a/n$,
\textit{Acta Arith.} 22 (1972/73) 339-352.
}

\bibitem{Vose}{
M. D. Vose,
Egyptian fractions,
\textit{Bull. London Math. Soc.} 17 (1985), no. 1, 21–24.
}

\begin{document}
\StartArtikel[Titel={Egyptian Fractions},
          AuteurA={Thomas F. Bloom},
          AdresA={Mathematical Institute\crlf 
          University of Oxford\crlf
          Woodstock Road\crlf 
          Oxford OX2 6GG, UK},
          EmailA={bloom@maths.ox.ac.uk},
          AuteurB={Christian Elsholtz},
          AdresB={Institut f\"{u}r Analysis und Zahlentheorie\crlf
                  Technische Universit\"{a}t Graz\crlf
                  Kopernikusgasse 24/II\crlf
                  A-8010 Graz, Austria},
          EmailB={elsholtz@math.tugraz.at},
	  ]
	 

\onderwerp{Introduction}
Scribes of Ancient Egypt had an unusual method of writing fractions: they were always\footnote{With the exceptions of $2/3$ and $3/4$.} expressed as a sum of distinct unit fractions $1/x_i$ -- for example, $\tfrac{2}{15}$ would instead be written as $\tfrac{1}{10}+\tfrac{1}{30}$. The Papyrus Rhind \cite{Rhind}, written over 3500 years ago, contains a table showing how to do this for the most frequently used fractions $2/n$. It is unclear precisely why the Egyptians wrote fractions in this form, which seems quite complicated for us, but the study of the equation  
\begin{equation}{\label{eq:general}}
    \frac{m}{n}=\frac{1}{x_1}+\cdots+\frac{1}{x_k}
\end{equation}
    with $m,n, x_1, \ldots, x_k$ positive integers, where in some cases $k$ is considered as fixed and in other cases $k$ may vary, 
has been very fruitful from a modern number theoretic point of view. Equations with integer or rational variables are called ``Diophantine equations". Another example of a Diophantine equation is the Pythagorean triple $x^2+y^2=z^2$, whose history goes back to Babylonian clay tablets, and is given a parametric solution in Euclid's ``Elements".

The study of solutions to \eqref{eq:general}, so-called `Egyptian fractions', has an important caveat when compared to the study of general Diophantine equations: we usually limit our attention to those solutions with $x_1 < x_2 < \cdots < x_k$ -- in particular, we do not allow equality between the $x_i$, and always count solutions as ordered tuples. For the study of Egyptian fractions this is no loss of generality: it has been observed by Takenouchi \cite{Takenouchi} that a fraction which can be written as a sum of $k$ unit fractions with repeated fractions can also be written as a sum of $k$ unit fractions with distinct fractions. To see this, one can replace any multiple occurrence of some $x_i$ by means of the two formulae below, depending on whether $x_i$ is odd or even:
$\frac{1}{2t}+\frac{1}{2t}=\frac{1}{t+1}+\frac{1}{t(t+1)}$ and $\frac{1}{2t+1}+\frac{1}{2t+1}=\frac{1}{t+1}+\frac{1}{(t+1)(2t+1)}$. Note that the sum of the denominators increases by this replacement. One possibly has to repeat this, but as the sum of the denominators increases, and as there is only a bounded number of ways to write a fraction as a sum of $k$ unit fractions (as we will show in Section~\ref{sec:count}), this is a bounded process which eventually stops with $k$ distinct unit fractions.

There are some even easier observations about Egyptian fractions: 
\begin{enumerate}
\item[a)] A fraction $m/n$ which has a solution of (\ref{eq:general}) with fixed $k$ also has a solution with every $k'\geq k$. 
\item[b)] If $m/n$ has a solution of the form (\ref{eq:general}) with fixed $k$, then $\frac{m}{nt}$ also has such a solution.
When $m$ and $k$ are fixed one therefore often concentrates on prime values $n$.
\end{enumerate}

The most basic question concerning Egyptian fractions is, given some rational $\frac{m}{n}\in (0,1]$, does \eqref{eq:general} always have a solution (for some $k\geq 1$)? The answer is yes, and in fact a greedy algorithm can be used to construct a solution: let $x_1\geq 1$ be the smallest positive integer such that $\frac{1}{x_1}\leq \frac{m}{n}$.
It follows that $\frac{m}{n}-\frac{1}{x_1}=\frac{m'}{n'}<\frac{1}{x_1}$ where $m'=mx_1-n<m$. While $\frac{m'}{n'}>0$ we may repeat this process, and since the numerator decreases at each stage this must terminate with a solution to \eqref{eq:general}. In particular, this greedy algorithm shows that $\frac{m}{n}$ can always be written as the sum of at most $m$ distinct unit fractions. 

For example, when applied to $4/17$ this algorithm produces the representation 
\[\frac{4}{17}=\frac{1}{5}+\frac{1}{29}+\frac{1}{1233}+\frac{1}{3039345}.\]
The solution found by the greedy algorithm is not necessarily the simplest one, 
in either the sizes of the denominators or the number of summands. For example, $4/17$ may also be written as
\[\frac{4}{17}=\frac{1}{6}+\frac{1}{17}+\frac{1}{102}.\]

The greedy approach is one of a variety of different algorithms for producing a solution to \eqref{eq:general}, each with its own advantages and disadvantages. For a survey on the rich literature concerning algorithmic aspects we refer the reader to \cite{Eppstein}.

In this survey we will consider more theoretical questions, and focus on a number of aspects that have received considerable interest over the last decades. We have endeavoured to give, where possible, some indication of the ideas and methods used, and hence have sacrificed some breadth for depth. The topics considered here do not cover the full range of results and open problems in this fascinating area, and the reader is encouraged to explore further in the surveys of Erd\H{o}s and Graham \cite{ErdosandGraham, Graham:2013}, the open problems collection of Guy \cite[Section D11]{Guy}, and the survey of Schinzel \cite{Schinzel}.

\onderwerp{The Erd\H{o}s-Straus conjecture}{\label{sec:erdos-straus}}
The greedy algorithm guarantees that $\frac{m}{n}$ can always be written as the sum of at most $m$ distinct unit fractions. As the example for $4/17$ above shows, however, this is not necessarily optimal, and one can ask, for any fraction $\frac{m}{n}$, what is the minimal $k$ such that \eqref{eq:general} has a solution.

\subsection{Sums of two unit fractions}
One certainly needs at least two unit fractions to represent $\frac{2}{n}$ (when $n$ is odd), and it is easy to check that one needs at least
three unit fractions to represent $\frac{3}{p}$, where $p$ is a prime such that $p\equiv 1 \bmod 3$. 

In general, Stewart \cite{Stewart:1964} has shown that a reduced fraction $\frac{m}{n}$ is a sum of two unit fractions if and only if there are two coprime divisors $n_1,n_2$ of $n$ such that $m$ divides $n_1+n_2$.

In particular, when applied to $m=4$, this implies that $\frac{4}{n}$ is the sum of two unit fractions for almost all $n$: Stewart's criterion implies that if $\frac{4}{n}$ is not the sum of two unit fractions then all prime factors of $n$ are $1 \bmod 4$, and the number of such integers $n \leq N$  is asymptotically 
$c \frac{N}{\sqrt{\log N}}=o(N)$, for some constant $c>0$.
In fact, for every $m$ a similar but more complicated analysis shows that the number of integers $n \leq N$ such that $\frac{m}{n}$ is not the sum of two unit fractions is asymptotically 
$c_m \frac{N(\log \log N)^{\beta_m}}{(\log N)^{\alpha_m}}= O(\frac{N}{\sqrt{\log N}})$, for some constants $c_m,\alpha_m,\beta_m$ -- see \cite{Elsholtz:1998} for details.

In particular, once we have fixed $m$, asymptotically almost all fractions with numerator $m$ are the sum of two unit fractions. The fact that we have fixed the numerator is vital here: for any fixed $k\geq 1$ the set of rationals which can be expressed as as sum of $k$ unit fractions is nowhere dense \cite{Sierpinski:1956}, except at $0$, so the picture changes considerably when not fixing the value $m$.

\subsection{Sums of three unit fractions}{\label{subsection:ES}}

When $m\leq 3$, the greedy algorithm produces a solution to \eqref{eq:general} with $k\leq 3$. When $m=4$ the greedy algorithm may require four unit fractions, as we have seen above. It is believed, however, that \emph{for all $n\geq 1$} this can be improved, and $\frac{4}{n}$ can be written as the sum of three unit fractions. This is perhaps the most notorious open problem concerning Egyptian fractions. 

\begin{conjecture}[Erd\H{o}s-Straus 1950]
For every $n\geq 2$ there exist positive integers $x,y,z$ such that
\begin{equation}\label{eq-es}
\frac{4}{n}=\frac{1}{x}+\frac{1}{y}+\frac{1}{z}.
\end{equation}
\end{conjecture}

This conjecture appeared in a 1950 paper of Erd\H{o}s \cite{Erdos:1950}, attributed to himself and Ernst G. Straus. It is unknown when the conjecture was originally made. Another reference 
is in a paper of Obl\'{a}th \cite{Obl1950}, submitted in 1948, who mentions it as a conjecture of Erd\H{o}s\footnote{Some webpages, such as Wikipedia, and a paper by Graham \cite{Graham:2013} seem to suggest that Obl\'{a}th made this conjecture independently. This is not the case, as Obl\'{a}th's paper clearly attributes it to Erd\H{o}s.},
and solved it for  small integers. When the second author asked Erd\H{o}s in 1996 how he came up with the conjecture
the answer was that this is the first interesting case.

Although $m=4$ is the first non-trivial case, it is believed that a similar phenomenon holds for any $m\geq 4$ (excluding finitely many exceptions). The analogous conjecture with numerator $5$ is due to Sierpi\'{n}ski
\cite{Sierpinski:1956}, and the general form was conjectured by Schinzel (also in \cite{Sierpinski:1956}).
\begin{conjecture}[Schinzel 1956]
For every $m\geq 4$ there exists a number $N_m$ such that, for every $n \geq N_m$,
there exist positive integers $x,y,z$ such that
\[\frac{m}{n}=\frac{1}{x}+\frac{1}{y}+\frac{1}{z}.\]
\end{conjecture}

We note that if \eqref{eq-es} is soluble for $n$ it is trivially also soluble for all multiples of $n$. In particular, in exploring the Erd\H{o}s-Straus conjecture it suffices to concentrate on prime values of $n$. 


The Erd\H{o}s-Straus conjecture has been computationally verified up to $10^{17}$ \cite{Sal}. It is possible to prove the Erd\H{o}s-Straus conjecture for many congruence classes via elementary means.

The key observation is that, for any $m\equiv 3\pmod{4}$, if there are integers $a,c,d$ such that $4acd-1=m$ and $n\equiv -a/c\pmod{m}$, then \eqref{eq-es} is solvable. Indeed, if $cn+a=(4acd-1)b$, then dividing by $abcdn$  shows that
\[\frac{4}{n} = \frac{1}{abd}+\frac{1}{acdn}+\frac{1}{bcdn}.\]
This observation allows one to verify the Erd\H{o}s-Straus conjecture for many congruence classes immediately. For example, modulo $47$ we could take $(a,c,d)=(1,6,2)$, so that any $n\equiv -8\pmod{47}$ satisfies the conjecture, and in fact we can write $4/47k-8$ as
\[{\textstyle \frac{1}{(6k-1)2}+\frac{1}{12(47k-8)}+\frac{1}{(6k-1)12 (47k-8)}.}\]
We could instead take $(a,c,d)=(1,2,6)$, thereby verifying the conjecture for $n\equiv -24\pmod{47}$, or $(a,c,d)=(2,3,2)$, verifying the conjecture for $n\equiv -32\pmod{47}$, and so on. Varying over all 18 ways of writing $12=acd$ we find 13 distinct congruence classes modulo $47$ for which the Erd\H{o}s-Straus conjecture is true. In general, the number of such `solved' congruence classes modulo $4t-1$ essentially depends on $d_3(t)$, the number of ways to write $t$ as a product of three positive integers, which on average grows like $\approx (\log t)^2$.

Using similar reasoning, one can quickly show that modulo 840 only the congruence classes $1,49,121,169,289,361$ (all squares!) are not generally solved (for details see \cite{Mor1969}), and modulo 120120 there is corresponding work by Terzi \cite{Ter1971}. 

If we take $a=c=1$ then the above shows that modulo all integers of the shape $4d-1$ the congruence class $-1$ is soluble. Any simple sieve approach (such as Brun's sieve) can then be used to prove that almost all integers $n \leq N$ satisfy the Erd\H{o}s-Straus conjecture.

A much more sophisticated version by Vaughan \cite{Vau1970} combines the above remark that the number of soluble congruence classes modulo primes is described by the divisor function $d_3$ with the large sieve and proves that there are at most $N\exp(-c (\log N)^{2/3})$ exceptions $n \in [1,N]$, for some absolute constant $c>0$. (To estimate the number of soluble congruence classes correctly, one needs both the Bombieri-Vinogradov and the Brun-Titchmarsh theorem, which guarantee that on average prime numbers are relatively well distributed in congruence classes.)

Arguing heuristically, the number of soluble congruence classes is so large that all sufficiently large integers should be covered by at least one, giving a compelling heuristic argument that the Erd\H{o}s-Straus conjecture holds, at least for all sufficiently large $n$. More precisely, one could argue as follows: since the average value of $d_3(m)$ is $\approx (\log m)^2$, for any prime $p\leq N$ congruent to $-1\pmod{4}$ we expect to find on average $\approx (\log p)^2$ many congruence classes modulo $p$ for which the Erd\H{o}s-Straus conjecture holds. Therefore, assuming independence between these classes for distinct primes, the `probability' that any $n$ fails the Erd\H{o}s-Straus conjecture is, using standard prime number estimates,
\begin{align*}
&\ll \prod_{\substack{p\leq n\\ p\equiv -1\pmod{4}}}\left(1-\frac{d_3((p+1)/4)}{p}\right)\\
&\approx e^{-\Omega((\log n)^2)}.
\end{align*}
(Here we use the Vinogradov notation $f \ll g$ to mean $f=O(g)$.) Since $\sum_{n\geq 1}e^{-c(\log n)^2}$ converges, this probabilistic heuristic suggests there are at most finitely many exceptions to the Erd\H{o}s-Straus conjecture. (Indeed, since this converges very rapidly, and the conjecture has already been confirmed up to $10^{17}$, this strongly suggests the conjecture holds for all $n$.)

Some congruence classes are easier than others for the Erd\H{o}s-Straus conjecture. For example, for $n \equiv 3 \bmod 4$ the fraction $\frac{4}{n}$ can even be represented by 2 unit fractions, $\frac{4}{4t+3}=\frac{1}{t+1}+\frac{1}{(t+1)(4t+3)}$,
and for $n\equiv 5 \bmod 8$ it can be represented by three unit fractions $\frac{4}{8t+5}=\frac{1}{2(t+1)}+\frac{1}{(t+1)(8t+5)}+\frac{1}{2(t+1)(8t+5)}$. It can be shown, however, that there are some congruence classes for which there is no such explicit formula. For example, Schinzel \cite{Schinzel:2000} proved this is the case for all quadratic residue congruence classes. In particular, there is no such formula for $n=4t+1$.

Although there are some congruence classes for which we cannot solve the Erd\H{o}s-Straus conjecture in such an explicit fashion, there are some congruence classes (as discussed above) for which we are able to easily verify the conjecture. To establish the conjecture for all primes $p$ (and hence all integers $n$) it therefore suffices to show that these `good' congruence classes cover all primes. (For example, although $61$ is of the form $4t+1$, for which we have no general solution, it is also of the form $7t-2$, for which we have the general solution $\tfrac{4}{7t-2}=\tfrac{1}{2t}+\tfrac{1}{2(7t-2)}+\tfrac{1}{t(7t-2)}$.)  This is a classical observation. It is less widely known that in fact this covering property is \emph{equivalent} to the Erd\H{o}s-Straus conjecture.

\begin{theorem}{\label{thm:equivalence}}
The Erd\H{o}s-Straus conjecture is equivalent to the statement that all primes are in at least one of the following congruence classes:
\[-a/c\pmod{4acd-1}\]
for some $a,c,d\geq 1$, or
\[-\frac{4c^2d+1}{k}\pmod{4cd}\]
for some $c,d,k\geq 1$ with $k\mid 4c^2d+1$.
\end{theorem}
(Similar statements can be found in Nakayama \cite{Nak1939}, Rosati \cite{Ros1954}, and Mordell \cite{Mor1969}.)
\begin{proof}
We first show that the covering statement is sufficient for the Erd\H{o}s-Straus conjecture to hold. The case $p\equiv -a/c\pmod{4acd-1}$ has already been discussed above. If $k\mid 4c^2d+1$ and $p\equiv -\frac{4c^2d+1}{k}\pmod{4cd}$ then there exists $a\geq 1$ such that $p=4acd-\frac{4c^2d+1}{k}$. Therefore $kp+1=4cd(ak-c)$ and an elementary rearrangement shows that
\[\frac{4}{p}=\frac{1}{ad(ak-c)}+\frac{1}{acd}+\frac{1}{(ak-c)cdp}.\]
Thus the Erd\H{o}s-Straus conjecture holds for all primes $p$ (and hence for all integers $n$) if these congruence classes cover all primes. 

We now argue that the covering statement is necessary. Suppose then that the Erd\H{o}s-Straus conjecture holds, let $p$ be some prime, and let $\frac{4}{p}=\frac{1}{x}+\frac{1}{y}+\frac{1}{z}$ with $x\leq y\leq z$, so that $4xyz = p(xy+yz+xz)$. Note that $x<p$, and hence $p\nmid x$. In fact one can show via an elementary argument, considering greatest common divisors, that there must exist integers $a,b,c,d\geq 1$ such that either
\begin{enumerate}
    \item $x=abd$, $y=acdp$, and $z=bcdp$, or
    \item $x=abd$, $y=acd$, and $z=bcdp$.
\end{enumerate}
In the first case we have $4abcdp^2=p(ap+cp^2+bp)$ whence $4abcd=a+b+cp$ and so $p\equiv -a/c\pmod{4acd-1}$, and in the second case we have $4abcdp=p(a+cp+bp)$, whence $4abcd=a+(b+c)p$. In particular $a$ divides $b+c$, say $b+c=ak$, and so 
\[kp+1=4cd(ak-c)=4akcd-4c^2d\]
whence $k\mid 4c^2d+1$ and $p\equiv -\frac{4c^2d+1}{k}\pmod{4cd}$. Thus $p$ belongs to at least one of the required congruence classes as claimed.
\end{proof}

Recently Bright and Loughran \cite{BL} have studied the Erd\H{o}s-Straus conjecture using tools from modern algebraic geometry, in particular showing that there is no Brauer-Manin obstruction to the existence of solutions to \eqref{eq-es}.

\subsection{Counting solutions}
For any existence problem there is a corresponding counting problem. Let $f(n)$ count how many representations $\frac{4}{n}$ has as the sum of three distinct unit fractions, so that the Erd\H{o}s-Straus conjecture is the statement that $f(n)>0$ for $n\geq 2$ (and is equivalent to the statement that $f(p)>0$ for all primes $p\geq 2$).

Elsholtz and Tao \cite{elsholtz-tao} have shown that
\[\sum_{p\leq N} f(p) = N(\log N)^{2+o(1)},\]
where the sum is over the primes $p\in[1,N]$. As the number of such primes is 
asymptotically $\frac{N}{\log N}$,
one can deduce that, on average, $f(p)= (\log p)^{3+o(1)}$. Elsholtz and Planitzer \cite{Elsholtz-Planitzer:2020} proved that for almost all 
integers $n \leq N$
\[f(n)\geq (\log n)^{\log 6+o(1)}.\] 
(Note that $\log 6\approx 1.79$.) Moreover for infinitely many $n$ the number of solutions is much larger, namely $f(n)\geq \exp ((\log 6+o(1) \frac{\log n}{\log \log n})$.
This is larger than one might expect at first sight,
and improves upon results of Elsholtz and Tao \cite{elsholtz-tao}. The crucial idea is to study those $n$ consisting of many small primes, 
where one can choose many divisors $d$ of $n$ such that $\frac{4}{n}-\frac{1}{d}$ still has very many solutions as a sum of two unit fractions.


A corresponding bound is also known even when we restrict to primes: all reduced congruence classes $e\pmod f$ contain primes that have many solutions, at least $\exp(c_f \frac{\log p}{\log \log p})$ for some constant $c_f>0$.

Furthermore, Elsholtz and Tao \cite{elsholtz-tao} established the pointwise upper bound 
$f(p) \leq p^{\frac{3}{5}+o(1)}$ for primes. This was generalized by Elsholtz and Planitzer \cite{Elsholtz-Planitzer:2020} to composite denominators: $f(n) \leq n^{\frac{3}{5}+o(1)}$.
It seems possible that a much better bound of $O_{\varepsilon}(n^{\varepsilon})$ holds.

\section{Bounding the number of fractions required}
As we have already observed, the greedy algorithm implies that any $\frac{m}{n}\in (0,1)$ is the sum of at most $m$ distinct unit fractions, and Schinzel's conjecture implies that in fact three unit fractions suffice, assuming $n$ is sufficiently large compared to $m$.

What if $m$ is large compared to $n$? The number of fractions required may grow with $n$: for example, an Egyptian fraction decomposition of $\frac{n-1}{n}$ requires $\Omega(\log\log n)$ distinct unit fractions. (This follows from the bounds given in Section~\ref{sec:count}.) The greedy algorithm shows that $n-1$ unit fractions always suffice, but much better bounds are known.

The proofs of all reasonable upper bounds follow a similar scheme. Suppose we have an increasing sequence $N_1<N_2<\cdots$ of positive integers such that any $1<n<N_k$ is the sum of at most $F(N_{k-1})$ distinct divisors of $N_k$, for some increasing function $F$. The relevance of such a sequence, as we will now show, is that it implies that any $\frac{m}{n}\in (0,1)$ can be written as the sum of at most $2F(n)$ distinct unit fractions. 

To see this, given $\frac{m}{n}\in (0,1)$ choose $k$ such that $N_{k-1}<n\leq N_k$, and let $\ell<N_k$ be such that $\frac{\ell}{N_k}\leq \frac{m}{n}<\frac{\ell+1}{N_k}$, whence $0\leq mN_k-n\ell<n\leq N_k$. Writing both $mN_k-n\ell$ and $\ell$ as the sum of at most $F(N_{k-1})\leq F(n)$ distinct divisors of $N_k$ and using the identity
\[\frac{m}{n}=\frac{mN_k-n\ell}{nN_k}+\frac{\ell}{N_k}\]
we obtain $\frac{m}{n}$ as the sum of at most $2F(n)$ distinct unit fractions as claimed.

A trivial example of such a sequence is to take $N_k=2^k$, which allows for $F(n)=O(\log n)$ (already a vast improvement over the $O(n)$ delivered by the greedy algorithm). Erd\H{o}s \cite{Erdos:1950} observed that a more efficient choice is to take $N_k$ to be the product of the first $k$ primes, which allows for $F(n)=O(\log n/\log\log n)$ instead. The best-known construction to date is due to Vose \cite{Vose}, who constructed an explicit sequence $N_k$ with $F(n)=O(\sqrt{\log n})$, yielding the following result. 
\begin{theorem}[Vose \cite{Vose}]
Any fraction $\frac{m}{n}\in(0,1)$ can be written as the sum of $O(\sqrt{\log n})$ many distinct unit fractions. 
\end{theorem}

Erd\H{o}s \cite{Erdos:1950} conjectured that in fact any $\frac{m}{n}\in(0,1)$ can be written as the sum of $O(\log\log n)$ many distinct unit fractions, which would be the best possible as the example of $\frac{n-1}{n}$ shows.

\onderwerp{Parametric solutions of $\frac{m}{n}=\frac{1}{x_1}+ \cdots + \frac{1}{x_k}$}
Let $E_{m,k}(N)$ count those $n\leq N$ such that $\frac{m}{n}$ is not the sum of $k$ unit fractions. As mentioned in section \ref{sec:erdos-straus}, Vaughan \cite{Vau1970} showed that $E_{4,3}(N)\leq N\exp(-c(\log N)^{2/3})$, for some positive $c$. This has been generalised by Elsholtz.
\begin{theorem}[Elsholtz \cite{Elsholtz:2001}]{\label{thm:elsholtz-thesis}}
Let $m>k\geq 3$ be positive integers. Then
\[E_{m,k}(N) \leq N\exp(-c_{m,k}(\log N)^{1-\frac{1}{2^{k-1}-1}})\]
for some absolute constant $c_{m,k}>0$.
\end{theorem}
Note that $m=4$ and $k=3$ recovers Vaughan's bound. Viola \cite{Viola}  previously established a similar bound with $\frac{1}{2^{k-1}-1}$ replaced by $\frac{1}{k-1}$, which Shen \cite{Shen} had improved to $\frac{1}{k}$.

A key idea in the proof of Theorem~\ref{thm:elsholtz-thesis}, and a useful tool in general for studying Egyptian fractions, is the realisation that solutions to the equation $\frac{m}{n}=\frac{1}{x_1}+\cdots+\frac{1}{x_k}$ can naturally be described by $2^k-1$ variables. The earliest reference containing this idea, which we are aware of, is Dedekind \cite{Dedekind:1931}. A detailed explanation of the general case is in \cite{Elsholtz:1998}, and shorter explanations are in \cite{Elsholtz:2001, elsholtz-tao, elsholtz-planitzer:2021}.


For concreteness, we will explain the idea first when $k=3$, and then when $k=4$. Suppose that $\frac{m}{n}=\frac{1}{x_1}+\frac{1}{x_2}+\frac{1}{x_3}$. Let $t_{123}=\mathrm{gcd}(x_1,x_2,x_3)$, $t_{12}=\mathrm{gcd}(x_1,x_2)/t_{123}$, and similarly for $t_{13}$ and $t_{23}$. Note that $t_{12},t_{13},t_{23}$ are all coprime in pairs. It is elementary to check that $t_{12}t_{13}t_{123}$ divides $x_1$, whence we can write $x_1=t_1t_{12}t_{13}t_{123}$, and similarly for $x_2$ and $x_3$. It follows from the definition that $t_1,t_2,t_3$ are also all coprime in pairs, and also that $\gcd(t_1,t_{23})=\gcd(t_2,t_{13})=\gcd(t_3,t_{12})=1$.

Furthermore, if we write $\frac{m}{n}$ as
\[{\textstyle \frac{1}{t_1 t_{12} t_{13} t_{123}}+
\frac{1}{t_2 t_{12} t_{23} t_{123}}+
\frac{1}{t_3 t_{13} t_{23} t_{123}}}\]
and multiply by common denominators we obtain
\[ m t_1 t_2 t_3 t_{12} t_{13} t_{23} t_{123}\]\[=n\left(t_1 t_2 t_{12}+
t_1 t_3 t_{13}+ t_2 t_3 t_{23} \right)
\]
and hence (assuming that $\gcd(m,n)=1$) each $t_i$ divides $n$. In particular, when $n$ is prime, this leaves only the possibilities that the $t_i$ are 1 or $p$. In other words, when $n$ is prime, one has $7-3=4$ free parameters. 

To further illustrate the idea, we now examine the case $k=4$. Every quadruple of four positive integers $(x_1,x_2,x_3,x_4)$ 
can be written as
\[
\begin{array}{ll}
 & x_1 =t_1 t_{12}t_{13} t_{14} t_{123} t_{124} t_{134} t_{1234}\\
     & x_2 =t_2 t_{12}t_{23} t_{24} t_{123} t_{124} t_{234} t_{1234}\\
     & x_3 =t_3 t_{13}t_{23} t_{34} t_{123} t_{134} t_{234} t_{1234}\\
     & x_4 =t_4 t_{14}t_{24} t_{34} t_{124} t_{134} t_{234} t_{1234},\\
\end{array}
\]
where 
\[t_{1234}=\gcd(x_1,x_2,x_3,x_4),\]
\[t_{123}=\gcd(x_1,x_2,x_3)/t_{1234}\]
(and similarly for other three indices), $t_{12}=\gcd(x_1,x_2)/(t_{1234} t_{123}t_{124})$, and so on. Crucially, any pair $t_I,t_J$ with $I\not\subseteq J$ and $J\not\subseteq I$, for example $t_{123},t_{124}$ or $t_{1},t_{23}$, must be coprime. Given a solution to $\frac{m}{n}=\frac{1}{x_1}+\frac{1}{x_2}+\frac{1}{x_3}+\frac{1}{x_4}$, multiplying by the common denominator as above, we see that when $n=p$ is prime, each $t_i$ must be a divisor of $p$. In particular, the family of solutions is naturally described by $2^4-1-4=11$ free parameters. In general, with the sum of $k$ unit fractions we have $2^k-k-1$ free parameters.

We may now attempt a generalisation of Vaughan's argument, finding a large collection of congruence classes that attempt to cover most $n\leq N$. Crucially, the number of degrees of freedom when constructing these congruence classes grows exponentially in $k$, leading ultimately (after a great deal of technical difficulty and further sieve estimates) to Theorem~\ref{thm:elsholtz-thesis}. (For comparison,  Viola \cite{Viola} and Shen \cite{Shen} used $k$ and $k+1$ free parameters respectively, resulting in their weaker bounds.)

For example, when $k=4$, one can now solve the classes
\[\left( m t_{12} t_{23} t_{24} t_{123} t_{124} 
t_{234} t_{1234} -1\right) r\]
\[{\scriptstyle - m t_{12} t_{123} t_{124} t_{1234} 
\left( t_{12} t_{23} t_{123} + t_{12} t_{24} t_{124} \right) }.\]
For comparison, in the case $k=4$ Viola made use of the fact that one can solve the classes
\[\left(m t_{123} t_{124} t_{1234} -1 \right) r -(t_{123}+t_{124}),\]
and Shen used
\[\left(m t_{123} t_{124} t_{234} t_{1234} -1 \right) r \]
\[-m t_{123} t_{124} t_{1234} (t_{123}+t_{124}). \]

When $k=3$ the number of soluble congruence classes was described by the $d_3$ function; now when $k=4$ it is described by the $d_7$ function. That is, modulo $q\equiv -1 \bmod m$ we obtain a soluble congruence class for each way of splitting $\frac{q+1}{m}$ into a product of $7$ factors.

There is a complication when $k>3$ (resulting in the final bound of Theorem~\ref{thm:elsholtz-thesis} containing a $2^{k-1}-1$ where one might expect $2^{k}-k-2$), since the soluble congruence classes are now described by sums of products, rather than a single product as in Vaughan's argument. In particular it is much more difficult to prove that different congruence classes thus obtained are actually distinct.

The reader may wonder if such parametric solutions could be used more generally in the study of Diophantine equations. The answer is yes, in principle, but in practice this general approach often simplifies to a much easier one,
as one can see with the Fermat equation $x^n+y^n=z^n$. Here one quickly sees that nothing is gained, as the variables can be assumed to be coprime.

\onderwerp{Counting solutions to $1=\frac{1}{x_1}+\cdots+\frac{1}{x_r}$}\label{sec:count}

We now turn our attention from the study of \eqref{eq:general} for fixed small $k$ and varying $m/n$, and consider the opposing situation in which we fix $m/n=1$ and study the solutions to
\begin{equation}\label{eq:one}
1=\frac{1}{x_1}+\cdots+\frac{1}{x_k}
\end{equation}
as $k$ varies. An easy inductive argument shows that the number of solutions to \eqref{eq:one} for fixed $k\geq 1$ is finite. We now present an elementary argument which gives an explicit upper bound as follows (for more details and slight improvements see \cite{Sandor:2003}).

Fix some solution to \eqref{eq:one}, ordered so that $x_1\leq \cdots \leq x_k$. For $0\leq m<k$ define $y_m\in\mathbb{Q}$ by
\begin{equation}\label{eq:ydef}
    1-\sum_{1\leq i\leq m}\frac{1}{x_i}=\frac{1}{y_m},
\end{equation}
so that $y_0=1$ and $y_{k-1}=x_k$, and for $1\leq m<k-1$ we have $y_m< x_{m+1}$. By definition, $\frac{1}{y_{m+1}}=\frac{1}{y_m}-\frac{1}{x_{m+1}}$, whence $y_{m+1}\leq y_{m}(y_{m}+1)$ for all $0\leq m< k-1$. In particular, if $(u_i)_{i\geq 1}$ is the sequence defined by $u_1=1$ and $u_{i+1}=u_i(u_i+1)$, then $y_m\leq u_{m+1}$ for $0\leq m<k$. Since the left-hand side of \eqref{eq:ydef} is at most $k/x_{m+1}$ we deduce that, for all $1\leq i\leq k$, we have $x_i\leq ku_i$. The sequence $u_n^{2^{-n}}$ is strictly increasing and tends to $c_0=\lim_{n\to \infty}u_n^{2^{-n}}=1.26408\cdots$, the so-called Vardi constant,
and hence $x_i\leq kc_0^{2^i}$. It follows immediately that the number of solutions to \eqref{eq:one} is at most
\[k^k\cdot c_0^{2^{k}-1}=c_0^{2^k(1+o(1))}.\]
In particular, this grows doubly exponentially with $k$. It is surprisingly difficult to come up with considerably better upper bounds. Small improvements are possible, by estimating the number of the first $k-2, k-3$ or $k-4$ values of $x_i$ trivially as above, but using non-trivial upper bounds for the number of ways of writing a fixed fraction as a sum of $2,3$ or $4$ unit fractions, respectively.
The best upper bound currently known is still doubly exponential in $k$, and is due to Elsholtz and Planitzer \cite{elsholtz-planitzer:2021}. They show (using $k-4$ and $4$ fractions) that the number of solutions to \eqref{eq:one} is at most
\[c_0^{2^k(\frac{1}{5}+o(1))}.\]
Although the upper bounds are all doubly exponential in $k$, it was an open problem for some time whether this was the true order of magnitude of the number of solutions to \eqref{eq:one}. That the number of solutions is indeed increasing with (essentially) doubly exponential growth was shown by Konyagin \cite{Kon2014}, who proved that the number of solutions to \eqref{eq:one} is at least
\[2^{c^{\frac{k}{\log k}}}\]
for some constant $c>0$. We will now sketch a variant of Konyagin's construction, which yields this lower bound for an increasing sequence of $k$, even if we further ask that all denominators are odd. For full details see \cite{Elsholtz:2015, Kon2014}.

We make crucial use of fractions of the shape $\frac{1}{2^n-1}$ with highly composite $n$. The important feature of these fractions is that the denominator $2^n-1$ has many divisors, and for every divisor $m\mid 2^{n}-1$ there is a decomposition of $1/2^n-1$ as
\[\frac{1}{2^n-1+m}+\frac{1}{2^n-1+(2^n-1)^2/m}.\]
This means that if we can find at least one representation of $1$ as the sum of $k-1$ distinct unit fractions, one of which is $\frac{1}{2^n-1}$, then there are at least $d(2^n-1)-2k$ many ways to write $1$ as the sum of $k$ distinct unit fractions. (Here the $-2k$ is to avoid counting representations with repeated denominators.)
We can write $1/2^n-1$ as both 
\[
 \frac{1}{2^{2 n}-1}+{\scriptstyle \left(
\frac{1}{2^{2 n}}
 +\frac{1}{2^{2 n} \left(2^{2 n}-1\right)}+\frac{1}{2^n+1}\right) }
\]
and
\[\frac{1}{2^{n+1}-1}+\]
\[{\scriptstyle \left(\frac{1}{(2^n-1)(2^{n+1}-1)}+
\frac{1}{2^{n+1}}+\frac{1}{2^{n+1}(2^{n+1}-1)}\right)}.\]
We may use these identities to find a solution to \eqref{eq:one} containing $\frac{1}{2^n-1}$ for any $n$ with $k=O(\log n)$. To do so one writes the number $n$ in binary and applies a combination of the two identities, to increase the exponent by one, or double it, respectively. For example
the number $29=16+8+4+1=11101_2$ can be reached
from $1$ as follows (in binary): $0,1,10,11,110,111,1110,11100,11101$. Thus, one can construct
a decomposition of $1$ into unit fractions including $\frac{1}{2^n-1}$ with $k=O(\log n)$ fractions.

We have therefore shown the existence of $\gg d(2^n-1)$ many solutions to \eqref{eq:one} in $k=O(\log n)$ many variables. If $n$ is the product of the first $r$ primes then $r\sim \log n/\log\log n\approx k/\log k$ and it can be shown (see \cite[Lemma 2.1]{Elsholtz:2015}) that $2^n-1$ has many prime factors: $\omega(2^n-1)\geq 2^r-6$. It follows that $d(2^n-1)\geq 2^{\omega(2^n-1)}\geq 2^{2^r-6}$. Combining these observations we have found $\geq 2^{c^\frac{k}{\log k}}$ many distinct solutions to \eqref{eq:one} with $k$ variables as required.





\onderwerp{Solutions to $1=\frac{1}{x_1}+\cdots+\frac{1}{x_k}$ with restricted denominators}

We now turn our attention from counting to an existence problem: what restrictions can we impose on the denominators in \eqref{eq:one} while still being guaranteed of finding a solution?

\subsection{Restrictions on the sizes of the denominators}
The first natural question is: does there exist a solution to \eqref{eq:one} such that the denominators are all large? Or all small? There are various ways to make this question precise. Henceforth we fix some $k\geq 1$ and order a solution to \eqref{eq:one} as $x_1<\cdots<x_k$.

We may first ask: how small can the largest denominator $x_k$ be? Erd\H{o}s observed that since $1\geq \sum_{0\leq j<k}\frac{1}{x_k-j}\sim \log(\frac{x_k}{x_k-k})$, we must have $x_k\geq (1+\tfrac{1}{e-1}+o(1))k$, and asked whether this was best possible. This was proved by Martin.

\begin{theorem}[Martin \cite{Mar}]
For any $k\geq 1$ there is a solution to \eqref{eq:one} such that 
\[1<x_1<\cdots<x_k \leq  (\tfrac{e}{e-1}+o_{k\to\infty}(1))k.\]
\end{theorem}

We now ask, on the other hand: how large can the smallest denominator $x_1$ be? By a similar argument to the above, we must have $x_1\leq (\tfrac{1}{e-1}+o(1))k$, and Erd\H{o}s asked whether this was best possible. This was proved, at least for infinitely many $k$, by Croot. In fact, Croot proves the following stronger result.

\begin{theorem}[Croot \cite{Cro2}]\label{th-croot1}
For any $N>1$ there exists some $k\geq 1$ and a solution to \eqref{eq:one} such that 
\[N<x_1<\cdots <x_k \leq (e+o_{N\to\infty}(1))N.\]
\end{theorem}

Notice that since the sum of all reciprocals in $(N,(e+o(1))N)$ is $1+o(1)$, we must have $k=(e-1+o(1))N$. It immediately follows that there are infinitely many $k$ and an accompanying solution to \eqref{eq:one} with $x_1\geq (\tfrac{1}{e-1}+o(1))k$, as required. (Note also that Croot's result implies Martin's for infinitely many $k$.) Both Croot's and Martin's results are more general than we have stated here, concerning decompositions of arbitrary rationals, and we refer to \cite{Mar,Cro2} for more details.

\subsection{Restrictions to arbitrary sets}
Of course, one may impose many more restrictions on the denominators than simple size bounds. Finding a solution to \eqref{eq:one} is a challenge even when the restrictions are very mild: for example, it is a non-trivial task to manually find a representation of $1$ as the sum of distinct unit fractions with all denominators odd (and $>1$).

There is one trivial obstruction that prevents a solution to \eqref{eq:one} within small sets: certainly no solution can exist with denominators in $A$ if $\sum_{n\in A}\frac{1}{n}<1$. Thus, for example, the set of integers in the interval $[N,2N]$ contains no solution to \eqref{eq:one}, no matter how large $N$ is. On the other hand, Theorem~\ref{th-croot1} shows that, for intervals, this trivial obstruction is the only one, and that the set of integers in $[N,(e+o(1))N]$ must contain a solution for all large $N$.

If we consider restricting the denominators to some (infinite) arithmetic progression, there are no longer any obvious obstructions. Indeed, we may find a solution to \eqref{eq:one} within any infinite arithmetic progression, as shown by van Albada and van Lint \cite{bAlLint}. Taken together, these positive results about finding solutions to \eqref{eq:one} in short intervals and arbitrary congruence classes can be seen as showing that there are no `local' obstructions to the existence of a solution to \eqref{eq:one}. 

It is therefore natural to conjecture, as Erd\H{o}s and Graham did in \cite{ErdosandGraham}, that the equation \eqref{eq:one} enjoys a Ramsey-type property: whenever the integers are finitely coloured, there exists a monochromatic solution to \eqref{eq:one}. This was proved by Croot.

\begin{theorem}[Croot \cite{Cro2003}]\label{th-croot}
For any $r\geq 1$, if the integers are arbitrarily coloured with $r$ many colours, then there must be a monochromatic solution to $1=\frac{1}{x_1}+\cdots+\frac{1}{x_k}$ with $1<x_1<\cdots<x_k$.
\end{theorem}

In fact, Croot proves a strong quantitative version: there exists a constant $C>1$ such that, for any $r\geq 1$, if the integers in $[2,C^r]$ are coloured with $r$ many colours then there must exist a monochromatic solution to \eqref{eq:one}. This exponential behaviour is the best possible, since using a greedy approach one can $r$-colour the integers in $[2,e^{(1+o(1))r})$ so that the sum of all reciprocals in each colour class will be less than one, and hence certainly there can be no solution to \eqref{eq:one}. It is an interesting and open problem to improve the value of Croot's constant $C$. (In \cite{Cro2003} Croot shows that $C=e^{167000}$ is sufficient for large $r$.)
For $r=2$ colours a Ph.D. thesis by Andreas Hipler \cite{Hipler:2002} proves that the interval $[2,208]$ has this property, and that $208$ is sharp. The proof involved nontrivial computer calculations, as one has to study many distinct colourings of certain crucial integers in this interval. 
If we colour the integers in $r$ many colours then at least one colour class has upper density\footnote{We recall that the upper density of $A\subset \mathbb{N}$ is defined as $\limsup_{N\to\infty}\frac{\lvert A\cap [1,N]\rvert}{N}$.} $\geq 1/r$. The following result, which was also conjectured by Erd\H{o}s and Graham (for example in \cite{Graham:2013}), is therefore a natural strengthening of Croot's colouring result. 

\begin{theorem}[Bloom 2021+ \cite{Blo2021}\footnote{The proof of Theorem~\ref{th-bloom} has now been formally computer-verified, using the Lean proof assisant, by Bloom and Mehta. The formal version of the proof can be accessed at \texttt{https://github.com/b-mehta/unit-fractions}.}]\label{th-bloom}
Any subset of the integers with positive upper density contains a solution to $1=\frac{1}{x_1}+\cdots+\frac{1}{x_k}$ with $1<x_1<\cdots<x_k$. 
\end{theorem}

The proof of Theorem~\ref{th-bloom} extends the method used by Croot to prove Theorem~\ref{th-croot}: a variant of the Hardy-Littlewood circle method, combined with an ingenious elementary combinatorial argument. Croot actually proves a density result for sufficiently `smooth' integers: in particular, that any positive density set of integers $A$, satisfying the additional constraint that all prime factors of $n\in A$ are $\leq n^{\frac{1}{4}-o(1)}$, contains a solution to \eqref{eq:one}. This immediately implies Theorem~\ref{th-croot}, since any $r$-colouring of all integers must also $r$-colour all such smooth integers.

This is not sufficient for an unrestricted density result since, for example, the set of all $n$ with a prime divisor $>n^{1/2}$ has positive density. The chief novelty of \cite{Blo2021} is that it improves the technical strength of Croot's argument so that the smoothness threshold of $\frac{1}{4}-o(1)$ is raised to $1-o(1)$. This suffices to prove an unrestricted density result, since almost all integers $n$ have no prime factors $>n^{1-o(1)}$. 

\subsection{A sketch of Croot's method}

We now present a sketch of Croot's method and the proofs of Theorems~\ref{th-croot} and \ref{th-bloom}. The actual proofs are quite technical, and the interested reader is referred to \cite{Cro2003} and \cite{Blo2021} for full details.

Suppose we are given some finite set of integers $A$ such that
\begin{enumerate}
    \item $A\subseteq [N,O(N)]$,
    \item $A$ is `$N^\theta$-smooth', in the sense that all prime factors of $n\in A$ satisfy $p\leq N^\theta$, and
    \item $\lvert A\rvert \gg N$.
\end{enumerate} 
Our goal is to find some $S\subseteq A$ such that $\sum_{n\in S}\frac{1}{n}=1$. We begin by noting that, combining properties (1) and (3) (with appropriate choices of constants), we may assume that $\sum_{n\in A}\frac{1}{n}\in (2-o(1),2)$ (after possibly discarding some elements of $A$). Since we then have $\sum_{n\in S}\frac{1}{n}<2$ for all $S\subseteq A$, it suffices to find some non-empty $S\subseteq A$ with $\sum_{n\in S}\frac{1}{n}\in \mathbb{Z}$. This trivial recasting has the advantage that it naturally leads to the possibility of using exponential sums. In fact, if we let $P=\mathrm{lcm}(A)$, then a simple exercise using orthogonality shows that the number of $S\subseteq A$ with $\sum_{n\in S}\frac{1}{n}\in\mathbb{Z}$ is exactly equal to
\begin{equation}\label{eq:count}
\frac{1}{P}\sum_{-\frac{P}{2}<r\leq \frac{P}{2}}\prod_{n\in A}(1+e(r/n)),
\end{equation}
where $e(x)=e^{2\pi ix}$. It suffices therefore to prove that \eqref{eq:count} is $\geq 2$. 

The $r=0$ summand contributes exactly $2^{\lvert A\rvert}/P$, and elementary number theory shows that the $N^\theta$-smoothness assumption implies $P=e^{O(N^\theta)}$, which is $2^{o(\lvert A\rvert)}$ provided $\theta\leq 1-o(1)$. This `main term' therefore contributes $2^{(1-o(1))\lvert A\rvert}$, which is much larger than we require. The contribution from other small $r$, those with $0<\lvert r\rvert\leq N/4$, say, is harder to calculate exactly, but an elementary calculation shows that the sign of $\Re\prod_{n\in A}(1+e(r/n))$ is $\cos(\pi r\sum_{n\in A}\frac{1}{n})\prod_{n\in A}\cos(\pi r/n)$. In particular, since $\sum_{n\in A}\frac{1}{n}\in (2-o(1),2)$, the contribution to \eqref{eq:count} from all $r$ with $0<\lvert r\rvert\leq N/4$ is non-negative, and hence can be discarded in our quest for a lower bound. 

In particular, to show that \eqref{eq:count} is $\geq 2$ as required, it suffices to show that whenever $N/4<\lvert r\rvert \leq P/2$, we have
\begin{equation}\label{eq:minor}
\left\lvert \prod_{n\in A}(1+e(r/n))\right\rvert =o(2^{\lvert A\rvert}/P).
\end{equation}
If $r\not\in [-N/8,N/8]\pmod{n}$, then, writing $r_n\in [-n/2,n/2]$ for the integer such that $r\equiv r_n\pmod{n}$, we have $e(r/n)=e(r_n/n)=e(\xi)$ for some $\xi\in(c,1/2)$ (with $c>0$ some absolute constant), and hence $\lvert 1+e(r/n)\rvert \leq 2/c'$ for some other absolute constant $c'>1$. Since $P\leq e^{O(N^{\theta})}$, to prove \eqref{eq:minor} it therefore suffices to find $\gg N^\theta$ many $n\in A$ such that $r\not\in [-N/8,N/8]\pmod{n}$. 

We may equivalently phrase this as saying that, if $I$ is the interval of width $N/4$ centred at $r$, then there are $\gg N^{\theta}$ many $n\in A$ that do not divide any $x\in I$.

At this point we have left any mention of exponential sums (and indeed of unit fractions) far behind: Croot's method transforms the problem into a purely combinatorial question concerning the interaction between intervals and multiples of elements in $A$. Finding a satisfactory answer to this question is a subtle and difficult affair, however. The overall idea is that if there is some interval $I$ which fails the requisite property, then one can construct an integer $x\in I$ with `too many' divisors.

To explain this combinatorial procedure, we introduce a new `measure' of sets of integers (depending on the fixed set $A$),
\[\mu(X)=\frac{1}{\log\log N}\sum_{p\in \mathcal{P}_A\cap X}\frac{1}{p},\]
where $\mathcal{P}_n$ is the set of primes dividing $n$, and $\mathcal{P}_X=\cup_{n\in X}\mathcal{P}_n$. It is also convenient to introduce the notation $A_q=\{ n\in A: q\mid n\}$, for any integer $q\geq 1$. The three crucial facts about this measure $\mu$ are that
\begin{enumerate}
    \item for any $X$, $\mu(X)\leq 1+o(1)$,
    \item for any $x$ and $q$, if 
\[\lvert \{ n\in A_q: n\mid x\}\rvert \geq (N/q)(\log N)^{-o(1)},\]then $\mu(\mathcal{P}_{x})\geq e^{-1}-o(1)$, and
    \item if $0<\lvert x_1-x_2\rvert \ll N$ then \[\mu(\mathcal{P}_{x_1}\cup\mathcal{P}_{x_2}) = \mu(\mathcal{P}_{x_1})+\mu(\mathcal{P}_{x_2})+o(1).\]
\end{enumerate}

Suppose now that there exists some interval $I$ of width $N/4$ such that all but $o(N^\theta)$ many $n\in A$ divide some $x(n)\in I$. We will argue that there must exist some $x\in I$ divisible by \emph{all} primes in $\mathcal{P}_A$, and hence divisible by $P=\mathrm{lcm}(A)$, which is an immediate contradiction if, as in our application, $I$ is an interval of width $N/4$ centred at some $r$ with $N/4< \lvert r\rvert \leq P/2$. 

Let $p\in \mathcal{P}_A$, and consider $A_p$ -- heuristically, we expect $\lvert A_p\rvert \approx \lvert A\rvert/p\gg N/p$, which we will assume henceforth. In particular, provided $N/p\gg N^\theta$, for almost all $n\in A_p$ there exists some $x(n)\in I$ divisible by $n$. In fact, after some divisor sleight of hand, one can ensure that there is some $x_p\in I$ which equals $x(n)$ for $\gg \lvert A_p\rvert(\log N)^{-o(1)}$ many $n\in A_p$. Therefore, by property (2) of $\mu$ above, we have $\mu(\mathcal{P}_{x_p})\geq e^{-1}-o(1)$. 

Therefore all $p\in \mathcal{P}_A$ have an associated $x_p\in I$ such that $\mathcal{P}_{x_p}$ has $\mu$-weight at least $e^{-1}>1/3$. Combining properties (1) and (3) of $\mu$ it follows that there cannot be three distinct such $x_p\in I$. If all $x_p$ are in fact identical, then we have found some $x\in I$ divisible by all $p\in\mathcal{P}$, and hence by $P$ as required (we are assuming for simplicity that $P$ is squarefree here). The only remaining possibility is that all $x_p$ are one of two distinct $x,y\in I$. In this case, we may find some large subset $A'\subseteq A$, such that all $n\in A'$ divide one of either $x$ or $y$, and we perform a similar iteration with $A'$ replacing $A$. 

The above analysis succeeds provided $N/p\gg N^\theta$ for all primes $p\in\mathcal{P}_A$. Since the $N^\theta$-smooth hypothesis implies $p\leq N^\theta$, this in turn is guaranteed provided $\theta <1/2$, and hence we have proved \eqref{eq:minor}. Therefore, assuming a smoothness threshold of $N^{1/2-o(1)}$, we have found some $S\subseteq A$ with $\sum_{n\in S}\frac{1}{n}=1$, and obtain Croot's smooth density result. (Croot works with $\theta<1/4$ rather than $<1/2$ to ease some of the many technical difficulties we have ignored in this sketch, but in principle his method works up to any $\theta<1/2$.)

To prove an unrestricted density result such as Theorem~\ref{th-bloom}, we need to raise this smoothness threshold to $\theta=1-o(1)$. The key idea of \cite{Blo2021} is to give up on the strong pointwise bound \eqref{eq:minor}, and instead only aim to prove an averaged version, which suffices for our purposes. We note that in fact the above combinatorial argument shows that, for any interval $I$ of width $O(N)$, there exists some $x\in I$ divisible by every $p\in\mathcal{P}_A$ with the property that all but $o(N/p)$ many $n\in A_p$ divide some $x\in I$. If there are many such $p$, then the ensuing bound on $\lvert \prod_{n\in A}(1+e(r/n))\rvert$ is quite weak (particularly when some of the primes can be as large as $N^{1-o(1)}$), but on the other hand, there are not many possible values of $r$ for which this can occur, since such $r$ must be $O(N)$-close to a fixed multiple of many primes simultaneously. On the other hand, if there are few such $p$, then we cannot restrict $r$, but can recover a pointwise bound comparable to \eqref{eq:minor} in strength. Trading off the two gains, we are able to establish a version of \eqref{eq:minor} that holds on average, which is sufficient.

\section{Applications}

We conclude by mentioning some intriguing applications of Egyptian fractions to other areas of pure mathematics.
\subsection{Finite group theory}
One simple yet surprising application is within finite group theory. It is natural to ask what limits the number of conjugacy classes imposes on the underlying group. By considering Egyptian fractions, Landau showed that there are only finitely many possibilities.
\begin{theorem}[Landau \cite{Lan1903}]
For any $k\geq 1$ there are only finitely many finite groups with exactly $k$ conjugacy classes. 
\end{theorem}

This is an elementary consequence of the bounds of Section~\ref{sec:count} on the size of the denominators of solutions to \eqref{eq:one}. It suffices to show that, if $G$ has exactly $k$ conjugacy classes, then $\lvert G\rvert \ll_k 1$.  Suppose that $G$ has conjugacy classes of sizes $m_1,\ldots,m_k$. Since these partition $G$, it follows that $\lvert G\rvert=m_1+\cdots+m_k$. On the other hand, each $m_i$ is a divisor of $\lvert G\rvert$, and hence, dividing by $\lvert G\rvert$, we have
\[1=\frac{1}{\lvert G\rvert/m_1}+\cdots+\frac{1}{\lvert G\rvert/m_k}.\]
It follows that, by the upper bound proved in Section~\ref{sec:count}, we have $\lvert G\rvert/m_i\leq kc_0^{2^k}$ for all $1\leq i\leq k$. Furthermore, since the identity is only conjugate to itself, for some $i$ we have $m_i=1$. It follows that $\lvert G\rvert \leq kc_0^{2^k}$ as required. (Note that this explicit bound implies, for example, that any finite group $G$ has $\gg \log\log \lvert G\rvert$ many distinct conjugacy classes.) 

Landau's lower bound estimate has been strengthened and extended in a number of ways; see, for example, \cite{BFM} for extensions and references.

\subsection{Polytopes}
A more recent connection is to discrete geometry. Recall that a polytope in $\mathbb{R}^d$ is the convex hull of some finite set of points in $\mathbb{Z}^d$. We say that $P$ is
\begin{itemize}
    \item integer-free if there are no points in $\mathbb{Z}^d$ in the interior of $P$,
    \item weakly maximal if it is integer-free and there is no integer-free polytope that strictly contains $P$, and
    \item strongly maximal if it is integer-free and there is no integer-free, closed, convex, $d$-dimensional, set of any sort that strictly contains $P$.
\end{itemize}

It is not obvious whether there can exist integer-free polytopes that are weakly maximal yet not strongly maximal. Indeed, there do not exist such polytopes in dimensions $\leq 3$. It was shown in \cite{NZ11} that there do exist such polytopes in all dimensions $\geq 4$. Recently, Averkov has shown that in fact there must exist \emph{many} such polytopes, by establishing a close link between counting polytopes which are weakly, not strongly, maximal, and Egyptian fraction decompositions of \eqref{eq:one}.
\begin{theorem}[Averkov \cite{Ave}]
Let $d\geq 6$. Up to affine equivalence, the number of weakly maximal integer-free polytopes in $\mathbb{R}^d$ that are not strongly maximal is at least the number of solutions to \eqref{eq:one} with $k=d-5$ variables.
\end{theorem}
In particular, using Konyagin's lower bound, we deduce that there are at least $2^{c^{d/\log d}}$ many weakly maximal polytopes that are not strongly maximal (for some constant $c>0$).

\subsection{Huffman codes}
The study of solutions to \eqref{eq:one} when the denominators $x_i$ are (not necessarily distinct) powers of a fixed integer $t$ has a close connection to coding theory. For our purposes, a $k$-code in an alphabet of size $t$ is simply a set of $k$ distinct $\{0,\ldots,t-1\}$-strings. For many practical applications, it is preferred that the code be prefix-free, that is, no string appears as an initial segment of any other (in particular this allows for instant decoding). The Kraft-McMillan inequality states that if a prefix-free $k$-code has string lengths $l_1,\ldots,l_k$ then $\sum \frac{1}{t^{l_i}}\leq 1$. 

Prefix-free codes with the average word length
as small as possible, also known as compact Huffman codes, therefore have string lengths that satisfy $\sum \frac{1}{t^{l_i}}=1$. In fact, given any solution to this equation, a corresponding compact Huffman code can be produced. For example, there are three essentially distinct ways (with the word lengths $l_i$ ordered by size) of writing $1$ as a sum of five reciprocal powers of $2$: 
$1
=\frac{1}{2}+ \frac{1}{4}+\frac{1}{8}+\frac{1}{16}+\frac{1}{16}
=\frac{1}{2}+ \frac{1}{8}+\frac{1}{8}+\frac{1}{8}+\frac{1}{8}
=\frac{1}{4}+ \frac{1}{4}+\frac{1}{4}+\frac{1}{8}+\frac{1}{8}.
$
These correspond to the compact Huffman codes $\{0, 10,110,1110,1111\}$, $\{0,100,101,110,111\}$, and $\{00,01,10,110,111\}$. This correspondence shows that the number of compact Huffman codes is (up to equivalence) the number of solutions to $\sum \frac{1}{t^{l_i}}=1$. This is also equivalent to counting other natural combinatorial objects, such as the number of nonequivalent canonical rooted trees, or the number of bounded degree sequences. For more information we refer to 
\cite{Elsholtz-Heuberger-Prodinger, Elsholtz-Heuberger-Krenn}.

\completepublications

\end{document}